\documentclass[11pt]{article}
\usepackage{amsfonts,amsmath,amsxtra}
\usepackage{latexsym}
\usepackage{amssymb}

\def\hybrid{\topmargin 0pt      \oddsidemargin 0pt
        \headheight 0pt \headsep 0pt
        \textwidth 16.5cm
        \textheight 23cm
        \marginparwidth 0.0in
        \parskip 5pt plus 1pt   \jot = 1.5ex}
\catcode`\@=11
\def\marginnote#1{}
\newcount\hour
\newcount\minute
\newtoks\amorpm
\hour=\time\divide\hour by60 \minute=\time{\multiply\hour by60
\global\advance\minute by-\hour}
\edef\standardtime{{\ifnum\hour<12 \global\amorpm={am}%
        \else\global\amorpm={pm}\advance\hour by-12 \fi
        \ifnum\hour=0 \hour=12 \fi
      \number\hour:\ifnum\minute<10 0\fi\number\minute\the\amorpm}}
\edef\militarytime{\number\hour:\ifnum\minute<10 0\fi\number\minute}

\def\draftlabel#1{{\@bsphack\if@filesw {\let\thepage\relax
   \xdef\@gtempa{\write\@auxout{\string
      \newlabel{#1}{{\@currentlabel}{\thepage}}}}}\@gtempa
   \if@nobreak \ifvmode\nobreak\fi\fi\fi\@esphack}
        \gdef\@eqnlabel{#1}}
\def\@eqnlabel{}
\def\@vacuum{}
\def\draftmarginnote#1{\marginpar{\raggedright\scriptsize\tt#1}}

\def\draft{\oddsidemargin -0.1truein
        \def\@oddfoot{\sl preliminary draft \hfil
        \rm\thepage\hfil\sl\today\quad\militarytime}
        \let\@evenfoot\@oddfoot \overfullrule 3pt
        \let\label=\draftlabel
        \let\marginnote=\draftmarginnote
\def\@eqnnum{{\rm (\theequation)}
\rlap{\kern\marginparsep\tt\@eqnlabel}%
\global\let\@eqnlabel\@vacuum}  }

\newfont{\Bbbb}{msbm7 scaled 1\@ptsize00}
\newcommand{\zs}{\raise-1pt\hbox{$\mbox{\Bbbb Z}$}}

scaled 1\@ptsize00

\font\sevenmsa=msam6 
\newfam\msafam
\textfont\msafam=\sevenmsa
\def\hexnumber@#1{\ifnum#1<10 \number#1\else
\ifnum#1=10 A\else\ifnum#1=11 B\else\ifnum#1=12 C\else \ifnum#1=13
D\else\ifnum#1=14 E\else\ifnum#1=15 F\fi\fi\fi\fi\fi\fi\fi}
\def\msa@{\hexnumber@\msafam}
\def\llcorner{\delimiter"4\msa@78\msa@78 }
\def\lrcorner{\delimiter"5\msa@79\msa@79 }
\mathchardef\blacktriangleright="3\msa@49
\mathchardef\blacktriangleleft="3\msa@4A \font\tenmsb=msbm10 scaled
1\@ptsize00
\newfam\msbfam
\textfont\msbfam=\tenmsb \scriptfont\msbfam=\tenmsb

\newdimen\Squaresize \Squaresize=14pt
\newdimen\Thickness \Thickness=0.5pt

\def\Square#1{\hbox{\vrule width \Thickness
   \vbox to \Squaresize{\hrule height \Thickness\vss
      \hbox to \Squaresize{\hss#1\hss}
   \vss\hrule height\Thickness}
\unskip\vrule width \Thickness} \kern-\Thickness}

\def\Vsquare#1{\vbox{\Square{$#1$}}\kern-\Thickness}

\def\numberbysection{\@addtoreset{equation}{section}
        \def\theequation{\thesection.\arabic{equation}}}
\numberbysection

\renewcommand{\theequation}{\thesection.\arabic{equation}}
\def\titlepage{\@restonecolfalse\if@twocolumn\@restonecoltrue\onecolumn
     \else \newpage \fi \thispagestyle{empty}\c@page\z@
        \def\thefootnote{\fnsymbol{footnote}} }

\def\endtitlepage{\if@restonecol\twocolumn \else  \fi
        \def\thefootnote{\arabic{footnote}}
        \setcounter{footnote}{0}}  
\relax

\hybrid
\makeatletter
\newdimen\normalarrayskip            
\newdimen\minarrayskip               
\normalarrayskip\baselineskip \minarrayskip\jot
\newif\ifold             \oldtrue            \def\new{\oldfalse}
\def\arraymode{\ifold\relax\else\displaystyle\fi}
\def\eqnumphantom{\phantom{(\theequation)}} 
\def\@arrayskip{\ifold\baselineskip\z@\lineskip\z@
     \else
     \baselineskip\minarrayskip\lineskip1\baselineskip\fi}

\def\@arrayclassz{\ifcase \@lastchclass \@acolampacol \or
\@ampacol \or \or \or \@addamp \or
   \@acolampacol \or \@firstampfalse \@acol \fi
\edef\@preamble{\@preamble
  \ifcase \@chnum
     \hfil$\relax\arraymode\@sharp$\hfil
     \or $\relax\arraymode\@sharp$\hfil
     \or \hfil$\relax\arraymode\@sharp$\fi}}

\def\@array[#1]#2{\setbox\@arstrutbox=\hbox{\vrule
     height\arraystretch \ht\strutbox
     depth\arraystretch \dp\strutbox
width\z@}\@mkpream{#2}\edef\@preamble{\halign \noexpand\@halignto
\bgroup \tabskip\z@ \@arstrut \@preamble \tabskip\z@ \cr}%
\let\@startpbox\@@startpbox \let\@endpbox\@@endpbox
    \if #1t\vtop \else \if#1b\vbox \else \vcenter \fi\fi
  \bgroup \let\par\relax
  \let\@sharp##\let\protect\relax
  \@arrayskip\@preamble}
\def\eqnarray{\stepcounter{equation}%
              \let\@currentlabel=\theequation
              \global\@eqnswtrue
              \global\@eqcnt\z@
              \tabskip\@centering              
              \let\\=\@eqncr
              $$%
            \halign to \displaywidth  \bgroup
             \eqnumphantom \@eqnsel
      \hskip\@centering                               
    $\displaystyle  \tabskip\z@ {##}$%
    &\global\@eqcnt\@ne \hskip 2\arraycolsep
         $ \displaystyle  \arraymode{##}$\hfil
    &\global\@eqcnt\tw@ \hskip 2\arraycolsep
         $\displaystyle\tabskip\z@{##}$\hfil
         \tabskip\@centering
    &{##}\tabskip\z@\cr}
\makeatother
\newcommand{\RR}{{\mathbb{R}}}

\newcommand{\ZZ}{{\mathbb{Z}}}

\def\IC{\mathbb{C}}

\def\IR{\mathbb{R}}
\def\IZ{\mathbb{Z}}

\def\CH {\mathcal{H}}

\def\CP {\mathcal{P}}
\def\CQ {\mathcal{Q}}

\def\ch{{\cal H}}

\def\a {{\alpha}}

\def\la{\lambda}

\def\e{\epsilon}
\def\pr {\partial}

\def\ch{{\rm ch}}

\def\Tr{{\rm Tr}}

\newtheorem{te}{Theorem}[section]

\newtheorem{prop}{Proposition}[section]
\newtheorem{cor}{Corollary}[section]
\newtheorem{lem}{Lemma}[section]
\newtheorem{ex}{Example}[section]

\newcommand{\proof}{\noindent {\it Proof}.\,\,}
\newcommand\bqa{\begin{eqnarray}}
\newcommand\eqa{\end{eqnarray}}
\def\be{\begin{eqnarray}\new\begin{array}{cc}}
\def\ee{\end{array}\end{eqnarray}}

\def\beq{\begin{equation}}
\def\eeq{\end{equation}}
\def\bse{\begin{subequations}}                
\def\ese{\end{subequations}}
\def\bp{\begin{pmatrix}}
\def\ep{\end{pmatrix}}

\def\i{\imath}

\def\stack#1#2{\raise0.7pt\hbox{$\mathrel{\mathop{#2}\limits^{#1}}$}}
\def\tr{\triangleright}
\def\tl{\triangleleft}
\def\sem{\mathsurround=0pt \raise1pt
\hbox{$\scriptscriptstyle>\!\!$}\:\!\!\tl}
\def\mes{\mathsurround=0pt \tr\!\:\!\raise0.8pt
\hbox{$\scriptscriptstyle\!\!<$}\,}
\def\]{\mathsurround=0pt ]\raise-2pt\hbox{$_\ast$}}

\def\<{\langle}
\def\>{\rangle}

\def\CQ{{\cal Q}}

\def\ch{{\cal H}}

\def\CH{\mathcal{H}}

\def\we{\raise-1pt\hbox{$\,\stackrel{\wedge}{,}\,$}}
\def\tr{{\rm tr}\,}
\def\Tr{{\rm Tr}\,}
\def\pr {\partial}

\newcounter{pac}[section]

\newcounter{pacc}[subsection]

0

\title{\bf On a classical limit of $q$-deformed Whittaker functions}
\begin{document}
\author{Anton Gerasimov, Dimitri Lebedev and Sergey Oblezin}
\date{}
\maketitle

\renewcommand{\abstractname}{}

\begin{abstract}
\noindent {\bf Abstract}.
We  provide a derivation of the Givental integral representation
of the classical $\mathfrak{gl}_{\ell+1}$-Whittaker function
as a limit $q\to 1$ of
the  $q$-deformed $\mathfrak{gl}_{\ell+1}$-Whittaker function
represented as a sum over the Gelfand-Zetlin patterns.

\end{abstract}
\vspace{5 mm}

\section*{Introduction}

The $q$-deformed $\mathfrak{gl}_{\ell+1}$-Whittaker functions can be
defined as eigenfunctions of the $q$-deformed
$\mathfrak{gl}_{\ell+1}$-Toda chains \cite{Ru}, \cite{Et}.
Among various eigenfunctions
there exists a special class of eigenfunctions with the support
in the positive Weyl chamber. By analogy with the classical case we call such
functions the class one $q$-deformed
$\mathfrak{gl}_{\ell+1}$-Whittaker functions. In \cite{GLO1} an
explicit representation of the class one
 $q$-deformed $\mathfrak{gl}_{\ell+1}$-Whittaker function as a sum over
the Gelfand-Zetlin patterns was proposed. This representation
has remarkable  integrality and positivity properties. Precisely
each term in the sum  is a positive integer multiplied by a
weight factor $q^{wt}$ and a character of the torus
$U_1^{\ell+1}$. This allows to represent  the $q$-deformed
$\mathfrak{gl}_{\ell+1}$-Whittaker function  as a
character of a $\IC^*\times U_{\ell+1}$-module (i.e. it allows a
categorification). The interpretation of the $q$-deformed
$\mathfrak{gl}_{\ell+1}$-Whittaker function as a character shall be
considered as a $q$-version of the Shintani-Casselman-Shalika
formula \cite{Sh}, \cite{CS}. Indeed in the limit $q\to 0$ the
$q$-deformed $\mathfrak{gl}_{\ell+1}$-Whittaker function can be
identified with the non-Archimedean Whittaker function and the
representation of the $q$-deformed Whittaker function as a character
reduces to the standard  Shintani-Casselman-Shalika formula for
non-Archimedean Whittaker function \cite{Sh}, \cite{CS}.

In the limit $q\to 1$ the $q$-Whittaker functions reproduces  the
classical Whittaker functions. It was pointed out in \cite{GLO1} that in this limit
an explicit sum type representation of the class one $q$-deformed
$\mathfrak{gl}_{\ell+1}$-Whittaker function turns into the Givental
integral representation for the class  one
$\mathfrak{gl}_{\ell+1}$-Whittaker function \cite{Gi} (see also
\cite{GKLO}). Thus the Givental integral representation shall be
considered as the Archimedean counterpart of the
Shintani-Casselman-Shalika formula (for more
details on this interpretation see \cite{GLO2}, \cite{GLO3}, \cite{G}).
In this  note we provide a precise description of  the $q\to 1$ limit
reducing the $q$-deformed $\mathfrak{gl}_{\ell+1}$-Whittaker
function to its classical analog and
explicitly demonstrate that the Givental integral representation
arises as a limit of the sum representation of the $q$-deformed
$\mathfrak{gl}_{\ell+1}$-Whittaker function. This result is
given by Theorem 3.1. The established relation between a sum over
the Gelfand-Zetlin patterns for $\mathfrak{gl}_{\ell+1}$
and the Givental integrals  for
$\mathfrak{gl}_{\ell+1}$  is a special case of a general relation
between the Gelfand-Zetlin patterns and the  Givental type integrals for
classical  series of Lie algebras \cite{GLO4}. This relation
elucidates the
identification of the Givental and the Gelfand-Zetlin graphs noticed
in \cite{GLO4}. The relation between the Gelfand-Zetlin and Givental
constructions described in this note should be also compared with the
duality type relation introduced in  \cite{GLO5}.
We are going to discuss the general form of the
relation between the Gelfand-Zetlin and the Givental
constructions for classical Lie algebras  elsewhere.

{\em Acknowledgments}: The authors are grateful to A. Borodin and
G. Olshanski for their interest in this work. The research was supported by  Grant
RFBR-09-01-93108-NCNIL-a. AG was  also partly supported by Science
Foundation Ireland grant. The research of SO was partially supported
by P. Deligne's 2004 Balzan Prize in Mathematics.

\section{ $q$-deformed  $\mathfrak{gl}_{\ell+1}$-Whittaker
function}

In this Section we recall the explicit construction
of the class one $q$-deformed $\mathfrak{gl}_{\ell+1}$-Whittaker
functions derived in  \cite{GLO1}.
Quantum  $q$-deformed $\mathfrak{gl}_{\ell+1}$-Toda chain (see
e.g. \cite{Ru}, \cite{Et})  is
defined by a set of $\ell+1$  mutually commuting functionally
independent quantum Hamiltonians $\CH_r^{\mathfrak{gl}_{\ell+1}}$,
$r=1,\ldots ,\ell+1$:
 \be\label{comm}
  \CH_r^{\mathfrak{gl}_{\ell+1}}(\underline{p}_{\ell+1})\,=\,\sum_{I_r}\,\bigl(
  \widetilde{X}_{i_1}^{1-\delta_{i_2-i_1,\,1}}\cdot\ldots\cdot
  \widetilde{X}_{i_{r-1}}^{1-\delta_{i_r-i_{r-1},\,1}}\cdot
  \widetilde{X}_{i_r}^{1-\delta_{i_{r+1}-i_r,\,1}}\bigr)
  T_{i_1}\cdot\ldots\cdot T_{i_r}\,,
 \ee
where $r=1,\ldots,\ell+1$ and $i_{r+1}=\ell+2$. The summation in
\eqref{comm} goes over all ordered subsets $I_r
=\{i_1<i_2<\cdots<i_r\} $ of $\{1,2,\cdots,\ell+1\}$. Here we  use
the  notations
$$
 T_if(\underline{p}_{\ell+1})=f(\underline{\widetilde{p}}_{\ell+1}),
\hspace{1.5cm}
 \widetilde{p}_{\ell+1,k}=p_{\ell+1,k}+\delta_{k,i},
$$
$$
 \widetilde{X}_i=1-q^{p_{\ell+1,i}-p_{\ell+1,i+1}+1},\hspace{0.5cm}
 i=1,\ldots,\ell,\hspace{1.5cm}\widetilde{X}_{\ell+1}=1.
$$
The corresponding eigenvalue problem can be written in the following
form:
 \be\label{eiglat}
  \ch_r^{\mathfrak{gl}_{\ell+1}}(\underline{p}_{\ell+1})
  \Psi^{\mathfrak{gl}_{\ell+1}}_{z_1,\ldots ,z_{\ell+1}}
  (\underline{p}_{\ell+1})\,=\,(\sum_{ I_r}\prod\limits_{i\in I_r} z_i)\,\,
  \Psi^{\mathfrak{gl}_{\ell+1}}_{z_1,\ldots ,z_{\ell+1}}
  (\underline{p}_{\ell+1}),
 \ee
and the   first nontrivial Hamiltonian  is given by
 \be\label{FirstHamiltonian}
 {\CH}_1^{\mathfrak{gl}_{\ell+1}}(\underline{p}_{\ell+1})\,=\,
 \sum\limits_{i=1}^{\ell}(1-q^{p_{\ell+1,i}-p_{\ell+1,i+1}+1})T_i\,+\,
 T_{\ell+1}.
\ee

One of the main results of \cite{GLO1} now  can be formulated as follows. Given
$\underline{p}_{\ell+1}= (p_{\ell+1,1},\ldots,p_{\ell+1,\ell+1})$
let us denote by $\CP^{(\ell+1)}(\underline{p}_{\ell+1})$ a  set of
collections of the integer parameters $p_{k,i}$, $k=1,\ldots,\ell$, $i=1,\ldots,k$
satisfying the Gelfand-Zetlin conditions $p_{k+1,i}\geq p_{k,i}\geq
p_{k+1,i+1}$. Let $\CP_{\ell+1,\ell}(\underline{p}_{\ell+1})$ be a
set of $\underline{p}_\ell=(p_{\ell,1},\ldots,p_{\ell,\ell})$,
$p_{\ell,i}\in \IZ$,  satisfying the conditions $p_{\ell+1,i}\geq
p_{\ell,i}\geq p_{\ell+1,i+1}$.

\begin{te}
A common solution of the eigenvalue problem (\ref{eiglat}) can be
written in the following form.  For $\underline{p}_{\ell+1}$ being
in the dominant  domain $p_{\ell+1,1}\geq\ldots\geq
p_{\ell+1,\ell+1}$, the solution is given by \be\label{main}
 \Psi^{\mathfrak{gl}_{\ell+1}}_{z_1,\ldots,z_{\ell+1}}
(\underline{p}_{\ell+1})\,=\,
\sum_{p_{k,i}\in\CP^{(\ell+1)}(\underline{p}_{\ell+1})}\,\,
\prod_{k=1}^{\ell+1} z_k^{\sum_i  p_{k,i}-\sum_i p_{k-1,i}}\,\,
\\ \times\frac{\prod\limits_{k=2}^{\ell}\prod\limits_{i=1}^{k-1}
(p_{k,i}-p_{k,i+1})_q!}
{\prod\limits_{k=1}^{\ell}\prod\limits_{i=1}^k
(p_{k+1,i}-p_{k,i})_q!\,\, (p_{k,i}-p_{k+1,i+1})_q!},\ee where  we
use the notation $(n)_q!=(1-q)...(1-q^n)$.  When
$\underline{p}_{\ell+1}$ is outside the dominant domain we set
$$
\Psi^{\mathfrak{gl}_{\ell+1}}_{z_1,\ldots,z_{\ell+1}}
(p_{\ell+1,1},\ldots,p_{\ell+1,\ell+1})\,=\,0.
$$
\end{te}

\begin{ex} Let $\mathfrak{g}=\mathfrak{gl}_{2}$,
$p_{2,1}:=p_1\in\ZZ$, $ p_{2,2}:=p_2\in\ZZ$ and $p_{1,1}:=p\in \ZZ$.
The function
$$
\Psi_{z_1,z_2}^{{\mathfrak gl}_2}(p_{1},p_{2}) =\sum_{p_{2}\leq
p\leq p_{1}}\frac{ z_1^{p} z_2^{p_{1}+p_{2}-p}}
{(p_1-p)_q!(p-p_2)_q!},\qquad  p_{1}\geq p_{2}\,,
$$
$$
\Psi_{z_1,z_2}^{{\mathfrak gl}_2}(p_{1},p_{2})=0, \qquad
p_{1}<p_{2}\,,
$$
is a  common eigenfunction of  mutually commuting Hamiltonians
$$
{\cal H}_1^{{\mathfrak gl}_2}\,=\,(1-q^{p_1-p_2+1})T_1+T_2,\qquad
{\cal H}_2^{{\mathfrak gl}_2}=T_1T_2.
$$
\end{ex}

The formula \eqref{main} can be easily rewritten in the
recursive form.
\begin{cor}
The following recursive relation holds
 \be\label{qtodarec}
  \Psi^{\mathfrak{gl}_{\ell+1}}_{z_1,\ldots,z_{\ell+1}}
  (\underline{p}_{\ell+1})\,=\,
  \sum_{\underline{p}_\ell\in\CP_{\ell+1,\ell}(\underline{p}_{\ell+1})}\,\,
  \Delta(\underline{p}_{\ell}) \,\,
  z_{\ell+1}^{\sum_ip_{\ell+1,i}-\sum_i p_{\ell,i}}\,\,
  Q_{\ell+1,\ell}(\underline{p}_{\ell+1},\underline{p}_{\ell}|q)
  \Psi^{\mathfrak{gl}_{\ell}}_{z_1,\ldots,z_{\ell}}(\underline{p}_{\ell}),
  \ee
where
 \be
  Q_{\ell+1,\ell}(\underline{p}_{\ell+1},\underline{p}_{\ell}|q)\,
  =\,\frac{1}{\prod\limits_{i=1}^{\ell} (p_{\ell+1,i}-p_{\ell,i})_q!\,\,
  (p_{\ell,i}-p_{\ell+1,i+1})_q!},\\
  \Delta(\underline{p}_{\ell})=
  \prod_{i=1}^{\ell-1}(p_{\ell,i}-p_{\ell,i+1})_q!\,\,\,.
 \ee
\end{cor}

The following representations of the class one $q$-deformed
$\mathfrak{gl}_{\ell+1}$-Whittaker function  are a consequence of the
positivity and integrality of the coefficients of the $q$-series expansions of
each term in the sum \eqref{main} (see \cite{GLO1} for details).

\begin{prop}\label{qCS} {\it (i)}. There exists a $\IC^*\times GL_{\ell+1}(\IC)$
  module $V$ such that
the common eigenfunction \eqref{main} of  the $q$-deformed Toda chain
allows the following representation for $p_{\ell+1,1}\geq
p_{\ell+1,2}\geq\ldots \geq p_{\ell+1,\ell+1}$: \be\label{intc}$$
\Psi^{\mathfrak{gl}_{\ell+1}}_{\underline{\lambda}}
(\underline{p}_{\ell+1})\,=\,\Tr_{V}\, q^{L_0}\prod_{i=1}^{\ell+1}
q^{\lambda_i H_i},\ee where  $H_i:=E_{i,i}$,
$i=1,\ldots,\ell+1$ are Cartan generators of
$\mathfrak{gl}_{\ell+1}={\rm Lie}(GL_{\ell+1})$ and $L_0$ is a
generator of ${\rm Lie}(\IC^*)$.

{\it (ii)}.  There exists a finite-dimensional
$\IC^*\times GL_{\ell+1}(\IC)$
  module $V_f$ such that the following representation holds
for $p_{\ell+1,1}\geq p_{\ell+1,2}\geq\ldots \geq p_{\ell+1,\ell+1}$:
\be\label{inrep1}
\widetilde{\Psi}^{\mathfrak{gl}_{\ell+1}}_{\underline{\lambda}}
(\underline{p}_{\ell+1})\,=\,\Delta(\underline{p}_{\ell+1})\,\,
\Psi^{\mathfrak{gl}_{\ell+1}}_{\underline{\lambda}}
(\underline{p}_{\ell+1})=\Tr_{V_f}\,q^{L_0}\prod_{i=1}^{\ell+1}
q^{\lambda_i\,H_i}. \ee  The module $V$
entering \eqref{intc} and the module $V_f$ entering \eqref{inrep1}
have a structure of modules under the action of (quantum)  affine
Lie algebras \cite{GLO1}.
\end{prop}

\section{Classical limit of $q$-deformed Toda chain}

In this Section we define  a limit $q\to 1$ of the $q$-deformed
$\mathfrak{gl}_{\ell+1}$-Toda chain reproducing the standard
$\mathfrak{gl}_{\ell+1}$-Toda chain. We provide an explicit
check that the first two generators of the ring of quantum
Hamiltonians of $\mathfrak{gl}_{\ell+1}$-Toda chain arise as a limit
of some combinations of the following quantum Hamiltonians of the
$q$-deformed Toda chain
$$
 \CH_1^{\mathfrak{gl}_{\ell+1}}(\underline{p}_{\ell+1}|q)\,=\,
 \sum\limits_{i=1}^{\ell}
 \Big(1\,-\,q^{p_{\ell+1,i}-p_{\ell+1,i+1}+1}\Big)\,T_i\,\,
 +\,\,T_{\ell+1},
$$
$$
\CH_{\ell+1}^{\mathfrak{gl}_{\ell+1}}(\underline{p}_{\ell+1}|q)\,
=\,T_1T_2\cdots T_{\ell+1}.
$$
Let us introduce the following parametrization:
 \be\label{param}
  q=e^{-\epsilon },\qquad
  p_{\ell+1,k}=(\ell+2-2k)m(\epsilon)+
  x_{\ell+1,k}\epsilon^{-1}.
 \ee
Here  $m(\epsilon)\in \IZ$ is given by
$$
 m(\epsilon)\, =\,-[\epsilon^{-1}\ln \,\epsilon],
$$
and  $[x]\in \IZ$ is the integer part of $x$.

\begin{prop} The following limiting  relations hold:
$$
 H^{\mathfrak{gl}_{\ell+1}}_1(\underline{x}_{\ell+1})\,
 =\,\lim_{\epsilon\to 0}\,\,\frac{1}{\epsilon}\,\,
 \Big[\CH_1^{\mathfrak{gl}_{\ell+1}}
(\underline{p}_{\ell+1}(x,\epsilon)|q(\epsilon))-(\ell+1)
 \Big]\,,
$$
 \be\nonumber
  H^{\mathfrak{gl}_{\ell+1}}_2(\underline{x}_{\ell+1})\,
  =-\,\lim_{\epsilon\to 0}\,\,\frac{1}{\epsilon^2}\,\,
  \Big[
\CH_1^{\mathfrak{gl}_{\ell+1}}(\underline{p}_{\ell+1}(x,\epsilon)|q(\epsilon))-
\CH_{\ell+1}^{\mathfrak{gl}_{\ell+1}}(\underline{p}_{\ell+1}(x,\epsilon)
|q(\epsilon))-\ell
\ee
$$
+\frac{1}{2}(\CH_{\ell+1}^{\mathfrak{gl}_{\ell+1}}(\underline{p}_{\ell+1}
(x,\epsilon)|q(\epsilon))-1)^2 \Big]\,
$$
where $H_i^{\mathfrak{gl}_{\ell+1}}$, $i=1,2$ are the
standard quantum Hamiltonians of the $\mathfrak{gl}_{\ell+1}$-Toda chain:
$$
 H^{\mathfrak{gl}_{\ell+1}}_1(\underline{x}_{\ell+1})=
\sum_{i=1}^{\ell+1}\frac{\pr}{\pr
 x_{\ell+1,i}},
$$
$$
 H^{\mathfrak{gl}_{\ell+1}}_2(\underline{x}_{\ell+1})
\,=\,-\frac{1}{2}\sum_{i=1}^{\ell+1}\frac{\pr^2}{\pr
  x_i^2}+\sum_{i=1}^{\ell} e^{x_{i+1}-x_{i}}.
$$
\end{prop}
\proof Using the fact that $\exp(\epsilon  [2(\epsilon )^{-1}\ln
(\epsilon )]=\epsilon^2(1+O(\epsilon^2/\ln \epsilon))$  we have
$$
 \CH_1^{\mathfrak{gl}_{\ell+1}}(\underline{p}_{\ell+1}|q)\,=\,(\ell+1)\,
 +\,\epsilon\sum_{i=1}^{\ell+1}\frac{\pr}{\pr x_{\ell+1,i}}
 +\epsilon^2\Big(\frac{1}{2}\sum_{i=1}^{\ell+1}\,
 \frac{\pr^2}{\pr x_{\ell+1,i}^2}\,
 -\,\sum_{k=1}^{\ell}e^{x_{\ell+1,k+1}-x_{\ell+1,k}}\Big)\,+\,O(\e^3),
$$
$$
 \CH_{\ell+1}^{\mathfrak{gl}_{\ell+1}}(\underline{p}_{\ell+1}|q)\,=\,
 1\,+\,\epsilon \sum_{i=1}^{\ell+1}\frac{\pr}{\pr x_{\ell+1,i}} +
 \frac{1}{2}\,\epsilon^2\,\,\sum_{i,j=1}^{\ell+1}
\frac{\pr^2}{\pr x_{\ell+1,i}\pr x_{\ell+1,j}}\,+\,O(\e^3).
$$
Now the limiting formulas  can be  straightforwardly verified.
We have
$$
\CH_1^{\mathfrak{gl}_{\ell+1}}(\underline{p}_{\ell+1}|q)-
\CH_{\ell+1}^{\mathfrak{gl}_{\ell+1}}(\underline{p}_{\ell+1}|q)-\ell=
 \,\epsilon^2\Big(-\frac{1}{2}\sum_{i\neq j}^{\ell+1}\,
 \frac{\pr^2}{\pr x_{\ell+1,i}\pr x_{\ell+1,j}}\,
 -\,\sum_{k=1}^{\ell}e^{x_{\ell+1,k+1}-x_{\ell+1,k}}\Big)\,+\,O(\e^3),
$$
$$
\frac{1}{2}(\CH_{\ell+1}^{\mathfrak{gl}_{\ell+1}}(\underline{p}_{\ell+1}|q)-1)^2=
 \,\frac{1}{2}\epsilon^2\Big(\sum_{i,j=1}^{\ell+1}\,
 \frac{\pr^2}{\pr x_{\ell+1,i}\pr x_{\ell+1,j}}\,\Big)\,+\,O(\e^3),
$$
and thus
$$
\CH_1^{\mathfrak{gl}_{\ell+1}}(\underline{p}_{\ell+1}|q)-
\CH_{\ell+1}^{\mathfrak{gl}_{\ell+1}}(\underline{p}_{\ell+1}|q)-\ell
+\frac{1}{2}(\CH_{\ell+1}^{\mathfrak{gl}_{\ell+1}}(\underline{p}_{\ell+1}|q)-1)^2
$$
$$
= \,\epsilon^2\Big(\frac{1}{2}\sum_{i=1}^{\ell+1}\,
 \frac{\pr^2}{\pr x_{\ell+1,i}^2}\,
 -\,\sum_{k=1}^{\ell}e^{x_{\ell+1,k+1}-x_{\ell+1,k}}\Big)\,+\,O(\e^3).
$$
$\Box$

It is easy to see that the eigenfunction problem \eqref{eiglat} is
transformed into the standard eigenfunction problem if we use the
following parametrization of the spectral
variables $z_i=e^{\i\,\epsilon\lambda_i},\,i=1,\ldots
,\ell+1$.

\section{Classical limit of class one Whittaker function}

In the limit $q\to 1$ defined in the previous Section the class one
solution \eqref{main} of  the $q$-deformed
$\mathfrak{gl}_{\ell+1}$-Toda chain should goes to the class one
solution of the classical $\mathfrak{gl}_{\ell+1}$-Toda chain. In
the classical setting an integral representation for  class one
$\mathfrak{gl}_{\ell+1}$-Whittaker function was constructed by
Givental \cite{Gi}, (see \cite{GKLO} for a choice of the contour
realizing class one condition) \bqa\label{giv}
\psi_{\lambda_1,\ldots,\lambda_{\ell+1}}^{\mathfrak{gl}_{\ell+1}}
(x_1,\ldots,x_{\ell+1})\,
=\,\int_{C}\prod\limits_{k=1}^{\ell}d\underline{x}_k\,\,
e^{\mathcal{F}^{\mathfrak{gl}_{\ell+1}}(x) }, \eqa and the function
$\mathcal{F}^{{\mathfrak{gl}}_{\ell+1}}(x)$ is given by
\bqa\label{intrep} \hspace{-1cm}
\mathcal{F}^{{\mathfrak{gl}}_{\ell+1}}(x)\,
=\,\imath\sum\limits_{n=1}^{\ell+1}
\lambda_n\Big(\sum\limits_{i=1}^n
x_{n,i}-\sum\limits_{i=1}^{n-1}x_{n-1,i}\Big)\,
-\,\sum\limits_{k=1}^{\ell}\sum\limits_{i=1}^k
\Big(e^{x_{k,i}-x_{k+1,i}}+e^{x_{k+1,\,i+1}-x_{k,i}}\Big)\,. \eqa
Here $C\subset N_+$ is a small  deformation of the subspace
$\RR^{\frac{(\ell+1)\ell}{2}}\subset \IC^{\frac{(\ell+1)\ell}{2}}\,
$ making the integral (\ref{intrep}) convergent. Besides, we use the
following notation:
$\underline{\lambda}=(\lambda_1,\ldots,\lambda_{\ell+1})$;
 $x_i:=x_{\ell+1,i},\,\,\,i=1,\ldots,\ell+1$.

The integral representation \eqref{giv} allows a recursive
presentation analogous to \eqref{qtodarec}
\be\label{itergiv}
\psi_{\lambda_1,\ldots,\lambda_{\ell+1}}^{\mathfrak{gl}_{\ell+1}}
  (x_1,\ldots,x_{\ell+1})\,
 =\,\int\limits_{\IR^{\ell}}d\underline{x}_{\ell}\,\,
  Q^{\mathfrak{gl}_{\ell+1}}_{\mathfrak{gl}_{\ell}}
  (\underline{x}_{\ell+1};\underline{x}_{\ell}; \lambda_{\ell+1})
  \psi_{\lambda_1,\ldots,\lambda_{\ell}}^{\mathfrak{gl}_{\ell}}
  (x_{\ell,\,1},\ldots,x_{\ell\ell}),
\ee
where
 \be\label{GivKernel}
  Q^{\mathfrak{gl}_{\ell+1}}_{\mathfrak{gl}_{\ell}}
  (\underline{x}_{\ell+1};\underline{x}_{\ell};\lambda_{\ell+1})
  =\,\exp\Big\{\i\lambda_{\ell+1}\Big(\sum_{i=1}^{\ell+1}
  x_{\ell+1,i}- \sum_{i=1}^{\ell}x_{\ell,\,i}\Big)\\
  -\,\sum_{i=1}^{\ell}\Big(e^{x_{\ell,\,i}-x_{\ell+1,\,i}}
  +e^{x_{\ell+1,\,i+1}-x_{\ell,\,i}}\Big)
  \Big\},
 \ee
 and  we assume $Q^{\mathfrak{gl}_{1}}_{\mathfrak{gl}_{0}}(x_{11};\,\la_1)
=e^{\imath \lambda_1x_{11}}$.

In the following we demonstrate that  in the previously defined
limit $q\to 1$ the class one $q$-deformed
$\mathfrak{gl}_{\ell+1}$-Whittaker function given by  the sum
\eqref{main} indeed turns into the classical class one
$\mathfrak{gl}_{\ell+1}$-Whittaker function given by the integral
representation \eqref{giv}. In particular iterative formula
\eqref{qtodarec} turns into \eqref{itergiv}. For this purpose we
need the following asymptotic of the $q$-factorials entering
\eqref{main}.
\begin{lem}\label{Lemma}
Let us introduce the following functions
$$
f_\alpha(y,\epsilon)=(y/\epsilon+\alpha m(\epsilon))_q!, \qquad
\alpha=1,2,
$$
where $m(\epsilon)=-[\epsilon^{-1}\ln\epsilon]$,
$q=e^{-\epsilon}$. Then for
$\epsilon\to +0$ the following expansions hold:
 \be\label{MAIN1}
  f_1(y,\epsilon)=e^{A(\epsilon)+e^{-y}+O(\epsilon)}\,;
 \ee
 \be\label{MAIN2}
  f_2(y,\epsilon)=e^{A(\epsilon)+O(\epsilon^{\alpha-1})},
 \ee
 where $A(\epsilon)=-\frac{\pi^2}{6}\frac{1}{\epsilon}-\frac{1}{2}\ln\frac{\e}{2\pi}\, $.
\end{lem}

\proof Taking into account the identity
$$
 \ln
 \prod_{n=1}^N(1-q^n)=\sum_{n=1}^N\ln(1-q^n)=-\sum_{n=1}^{N}\sum_{r=1}^{+\infty}
 \frac{1}{r}q^{nr}=
 -\sum_{r=1}^{+\infty}\frac{q^r}{r}\left(\frac{1-q^{Nr}}{1-q^r}\right),
$$
and using the substitution  $q=e^{-\e}$, $N=\e^{-1}y+\alpha m(\epsilon)$ we
obtain
$$
 \ln\,f_\alpha(y,\epsilon)=
 -\sum_{r=1}^{+\infty}\frac{e^{-r\epsilon}}{r}\left(\frac{1-e^{-\alpha
       r \epsilon m(\epsilon)} e^{-ry}}
 {1-e^{-r\epsilon}}\right).
$$
Now expanding the denominator over small $\epsilon$ we have
$$
\ln\,f_\alpha(y,\epsilon)=
-\sum_{r=1}^{+\infty}\frac{e^{-r\epsilon}}{r}\left(\frac{1-\epsilon^{\alpha
      r}e^{-ry}}
{1-e^{-r\epsilon}}\right)+\cdots=
-\sum_{r=1}^{+\infty}\frac{e^{-r\epsilon}}{r^2\epsilon}\left(\frac{1-\epsilon^{\alpha
      r}e^{-ry}}
{1-\frac{1}{2}r\epsilon+\frac{1}{3!}r^2\epsilon^2+\cdots}\right)+\cdots,
$$
and for the derivative we obtain
$$
 \pr_y\ln\,f_{\alpha}(y,\epsilon)\,
 =\,-\sum_{r=1}^{+\infty}\frac{1}{r\epsilon}
 \left(\frac{\epsilon^{\alpha r}e^{-ry-r\epsilon}}
 {1-\frac{1}{2}r\epsilon+\frac{1}{3!}r^2\epsilon^2+\cdots}\right)+\cdots=
\sum_{k=-1}^{+\infty} c_k\,I_{\alpha,k}(y,\epsilon),
$$
where
$$
I_{\alpha,k}(y,\epsilon)\,=\,\sum_{r=1}^{+\infty} \epsilon^{k+\alpha r} r^k e^{-yr}\,=\,
\epsilon^k\sum_{r=1}^{+\infty} t^{r} r^k ,\qquad
t=e^{-y}\epsilon^{\alpha},
$$
and $c_{-1}=-1$. Let us separately analyze the term $I_{\alpha,-1}$ and the
other terms $I_{\alpha,k\geq 0}$.
We have $$
 I_{\alpha,k\geq 0}(y,\epsilon)\,
 =\,\epsilon^k\left(t\frac{\pr}{\pr t}\right)^k\frac{1}{1-t}\,
 =\,\epsilon^k\frac{\pr^k}{\pr y^k}\frac{1}{1-\epsilon^{\alpha}e^{-y}}\,,
$$
and thus
$$
I_{\alpha,k\geq0}\,=\,\epsilon^{k+\alpha}e^{-y}+\cdots,\qquad
\alpha=1,2.
$$
Now consider the case of $k=-1$
$$
c_{-1}I_{\alpha,-1}(y,\epsilon)=-\frac{1}{\epsilon}\,\sum_{r=1}^{+\infty}
\frac{t^{r}}{r}=-\frac{1}{\epsilon}\ln(1-t)=-\frac{1}{\epsilon}\ln(1-\epsilon^{\alpha}
e^{-y}))=\epsilon^{\alpha-1}e^{-y}+\cdots ,\qquad t=e^{-y}
\epsilon^{\alpha}.
$$
This gives \eqref{MAIN1}, \eqref{MAIN2}  with an  unknown $A(\epsilon)$. To
calculate $A(\epsilon)$ we take $e^{-y}=0$ and notice that the
resulting function does not depend on $\alpha$. Thus we
should calculate the asymptotic of the following function:
$$
 \ln\,f_{\a}(y,\epsilon)|_{e^y=0}=
 -\sum_{r=1}^{+\infty}\frac{e^{-r\epsilon}}{r}\left(\frac{1-\epsilon^{r\a}e^{-ry}}
 {1-e^{-r\epsilon}}\right)\Big|_{e^y=0} =-\sum_{r=1}^{+\infty}\frac{1}{r}
 \left(\frac{e^{-r\epsilon}}{1-e^{-r\epsilon}}\right)=
$$
$$
 =-\sum_{n=1}^{+\infty}\sum_{r=1}^{+\infty}
 \frac{1}{r}e^{-nr\epsilon}=
 \ln\prod_{n=1}^{+\infty}(1-e^{-n\epsilon}).
$$
It can be easily done using the modular properties
$$
 \eta(-\tau^{-1})\,=\,\sqrt{-\i\tau}\,\,\eta(\tau),
$$
of the Dedekind eta function
$$
 \eta(\tau)\,
 =\,e^{\frac{\imath\pi\tau}{12}}\,\,\prod_{n=1}^{\infty}
 (1-e^{2\pi\imath n\tau}).
$$
Namely, taking $\tau=\frac{\i\epsilon}{2\pi}$ we have
$$
 f_\alpha\bigl(y;\,\e\bigr)\Big|_{e^{-y}=0}\,
 =\,\sqrt{2\pi\e^{-1}}\,e^{-\frac{\pi^2}{6}\e^{-1}}
  \prod_{n=1}^{\infty}\bigl(1-e^{-\e^{-1}(2\pi)^2n}\bigr)\,.
$$
This allows to infer the following  result for the leading
coefficients in the asymptotic expansion of $\ln
f_{\alpha}(y,\e\bigr)\Big|_{e^{-y}=0}$:
 \be
   A(\epsilon)\,
  =\,-\frac{1}{2}\ln\frac{\e}{2\pi}-\frac{\pi^2}{6}\e^{-1}\,,
 \hspace{1.5cm}
  \e\longrightarrow+0\,.
 \ee
This completes the proof of Lemma. $\Box$

\begin{te} Let us use the following parametrization
\be\label{param1}
  q=e^{-\epsilon},\qquad
 p_{\ell+1,k}= (\ell+2-2k)m(\epsilon)+
  \epsilon^{-1}x_{\ell+1,k}\,,\qquad z_k=e^{\i\,\epsilon\lambda_k}\,,
 \ee
where $k=1,\ldots,\ell+1$, $m(\e)\,=\,-\bigl[\e^{-1}\ln\e\bigr]$.
 The  integral representation \eqref{giv} of the classical
$\mathfrak{gl}_{\ell+1}$-Whittaker function  is given by the following
limit of the $q$-deformed class one
$\mathfrak{gl}_{\ell+1}$-Whittaker function represented as a sum \eqref{main}
 \be\label{mainlim}
  \psi_{\lambda_1,\ldots,\lambda_{\ell+1}}^{\mathfrak{gl}_{\ell+1}}
  (x_1,\ldots,x_{\ell+1})\,
  =\lim_{\epsilon\to +0}\Big[\epsilon^{\frac{\ell(\ell+1)}{2}}\,
  e^{\frac{\ell(\ell+3)}{2}A(\e)}\,\,
  \Psi^{\mathfrak{gl}_{\ell+1}}_{z_1,\ldots ,z_{\ell+1}}
  (\underline{p}_{\ell+1})\Big]\,,
 \ee
where  
$A(\epsilon)=-\frac{\pi^2}{6}\frac{1}{\epsilon}-\frac{1}{2}\ln\frac{\e}{2\pi}$ and
 $x_i=x_{\ell+1,i},\,i=1,\ldots,\ell+1$.
\end{te}

\proof We prove \eqref{mainlim} by relating
the recursive relation \eqref{qtodarec}
\be\label{qtodarec1}
  \Psi^{\mathfrak{gl}_{\ell+1}}_{z_1,\ldots,z_{\ell+1}}
  (\underline{p}_{\ell+1})\,=\,
  \sum_{\underline{p}_\ell\in\CP_{\ell+1,\ell}(\underline{p}_{\ell+1})}
  \!\!\!
  \Delta(\underline{p}_{\ell}) \,\,
  z_{\ell+1}^{\sum\limits_{i=1}^{\ell+1}p_{\ell+1,i}
  -\sum\limits_{j=1}^{\ell}p_{\ell,j}}
  Q_{\ell+1,\ell}(\underline{p}_{\ell+1},\underline{p}_{\ell}|q)\,\,
  \Psi^{\mathfrak{gl}_{\ell}}_{z_1,\ldots,z_{\ell}}(\underline{p}_{\ell}),
 \nonumber
 \ee
where
 \be
  Q_{\ell+1,\ell}(\underline{p}_{\ell+1},\underline{p}_{\ell}|q)\,
  =\,\frac{1}{\prod\limits_{i=1}^{\ell} (p_{\ell+1,i}-p_{\ell,i})_q!\,\,
  (p_{\ell,i}-p_{\ell+1,i+1})_q!},\\
  \Delta(\underline{p}_{\ell})=
  \prod_{i=1}^{\ell-1}(p_{\ell,i}-p_{\ell,i+1})_q!\,,
 \ee
with the recursive relation
\eqref{itergiv} for the classical Whittaker function.

Let us introduce the following parametrization of the elements
of the Gelfand-Zetlin patterns
$\underline{p}_{\ell}\in\CP_{\ell+1,\,\ell}(\underline{p}_{\ell+1})$:
 \be\label{ScalingVariables}
  p_{\ell,\,k}\,=\,\e^{-1}x_{\ell,\,k}\,
  +\,a_k m(\e)\,,
 \hspace{1cm}
  m(\e)\,=\,-\bigl[\e^{-1}\ln\e\bigr],
 \ee
where $a_k$ are some constants.
The Gelfand-Zetlin conditions on weights
$\underline{p}_{\ell+1}$ reads as follows:
$$
 p_{\ell+1,\,k}\geq p_{\ell,\,k}\geq p_{\ell+1,\,k+1}\,,
\hspace{1.5cm}
 k=1,\ldots,\ell\,,
$$
and they lead to
 \bqa \label{uniq}
  \e^{-1}x_{\ell+1,\,k}+(\ell+2-2k)m(\e)\,
  \geq\,
  \e^{-1}x_{\ell,\,k}+a_km(\e)
 \geq\,
  \e^{-1}x_{\ell+1,\,k+1}+(\ell-2k)m(\e)\,.
 \eqa
The requirement that  the limit $\epsilon\to +0$ preserves the
conditions (\ref{uniq})  implies the following restrictions on the
parameters $a_k$:
\be
\ell-2k+2\,\,>\,\,
a_k\,\,>\,\,
\ell-2k\,,
\hspace{1.5cm}
k=1,\ldots,\ell\,.
\ee
Since $\underline{p}_{\ell}=(p_{\ell,1},\ldots ,p_{\ell,\ell})\in
\IZ^{\ell}$  the only consistent choice in the limit $\epsilon\to + 0$
is $a_k=\ell+1-2k$, $k=1,\ldots,\ell$.
Although the variables $p_{\ell,\,k}$ are restricted to be in positive
Weyl chamber i.e. $p_{\ell,k}\geq p_{\ell,k+1}$,
in the limit $\epsilon\to +0$ the variables $x_{\ell,\,k}$
have no such restrictions. This follows from
a simple observation that  the limit $\epsilon\to +0$
the $\frac{a}{\epsilon}-b[\epsilon^{-1}\ln\epsilon]\to
+\infty/-\infty$ depends only on the sign of non-zero coefficient $b$.
Thus we have
\be\label{ScalingVariables1}
 p_{\ell,\,k}\,=\,\e^{-1}x_{\ell,\,k}\,
  +\,(\ell+1-2k) m(\e)\,.
\ee
Now using Lemma \ref{Lemma}  it is easy to obtain the following
limiting formulas:
 \be
  \lim_{\e\to+0}e^{2\ell A(\e)}\CQ_{\ell+1,\ell}
  (\underline{p}_{\ell+1},\underline{p}_{\ell}|q)\,
  =\,\lim_{\e\to+0}
  \frac{e^{2\ell A(\e)}}{\prod\limits_{i=1}^{\ell}
  f_1(x_{\ell+1,\,i}-x_{\ell,\,i},\epsilon)\,\,
  f_1(x_{\ell,\,i}-x_{\ell+1,\,i+1},\epsilon)}\\
  =\,Q^{\mathfrak{gl}_{\ell+1}}_{\mathfrak{gl}_{\ell}}
  \bigl(\underline{x}_{\ell+1};\underline{x}_{\ell};\,
  \lambda_{\ell+1}\bigr)\Big|_{\lambda_{\ell+1}=0}\,,
 \ee
 \be
  \lim_{\e\to+0}e^{(1-\ell)A(\e)}\Delta(\underline{p}_{\ell})\,
  =\,\lim_{\e\to+0}e^{(1-\ell)A(\e)}
  \prod_{i=1}^{\ell-1}f_2(x_{\ell,i}-x_{\ell,i+1},\epsilon)\,\,=\,\,1\,,
 \ee
where $Q^{\mathfrak{gl}_{\ell+1}}_{\mathfrak{gl}_{\ell}}
  \bigl(\underline{x}_{\ell+1};\underline{x}_{\ell};\,
  \lambda_{\ell+1}\bigr)$ is given by \eqref{GivKernel}.
This implies the following identity:
 \be\nonumber
  \lim_{\e\to+0}\Big\{\e^{\ell}\!\!\!\!
  \sum_{\underline{p}_\ell\in\CP_{\ell+1,\ell}(\underline{p}_{\ell+1})}
  \!\!\!
  z_{\ell+1}^{\sum\limits_{i=1}^{\ell+1}p_{\ell+1,i}
  -\sum\limits_{j=1}^{\ell}p_{\ell,j}}\,\,
  \Big[e^{(\ell+1)A(\e)}Q_{\ell+1,\,\ell}
  \bigl(\underline{p}_{\ell+1};\underline{p}_{\ell}\bigr|q\bigr)\,
\Delta(\underline{p}_{\ell})\Big]\\
\times \,\,\,\e^{\frac{\ell(\ell-1)}{2}}e^{\frac{(\ell-1)(\ell+2)}{2}A(\epsilon)}
  \Psi^{\mathfrak{gl}_{\ell}}_{z_1,\ldots,z_{\ell}}
  (\underline{p}_{\ell})\Big\}
 \ee
 \be\label{LimitToIntegral}
  =\,\int\limits_{\IR^{\ell}}d\underline{x}_{\ell}\,\,
  \exp\Big\{\i\la_{\ell+1}\bigl(\sum\limits_{i=1}^{\ell+1}x_{\ell+1,i}
  -\sum\limits_{j=1}^{\ell}x_{\ell,j}\bigr)\Big\}\\
  \times\lim_{\e\to+0}\Big[e^{(\ell+1)A(\e)}Q_{\ell+1,\,\ell}
  \bigl(\underline{p}_{\ell+1}(\underline{x}_{\ell+1},\,\e);\,
  \underline{p}_{\ell}(\underline{x}_{\ell},\,\e)\,\bigr|\,q(\e)\bigr)\,\,
  \Delta\bigl(\underline{x}_{\ell}(\underline{x}_{\ell},\,\e)\bigr)\Big]\\
  \times
  \lim_{\e\to+0}\Big[
  \e^{\frac{\ell(\ell-1)}{2}}e^{\frac{(\ell-1)(\ell+2)}{2}A(\epsilon)}
  \Psi^{\mathfrak{gl}_{\ell}}_{z_1,\ldots,z_{\ell}}
  \bigl(\underline{p}_{\ell}(\underline{x},\,\e)\bigr)\Big]\,.
 \ee
Thus we recover the recursive relations \eqref{itergiv} for the Givental integrals
directly leading to the integral representation \eqref{giv} for the classical
$\mathfrak{gl}_{\ell+1}$-Whittaker function. Using
\eqref{LimitToIntegral} iteratively over $\ell$ we obtain
\eqref{mainlim}. $\Box$

\begin{ex} For $\ell=1$ we have
$$
\Psi^{\mathfrak{gl}_2}_{z_1,z_2}(p_{2,1},p_{2,2})
=\sum_{p_{2,2}\leq p_{1,1}\leq p_{2,1}}\frac{ z_1^{p_{1,1}}
z_2^{p_{2,1}+p_{2,2}-p_{1,1}}}
{(p_{1,1}-p_{2,2})_q!(p_{2,1}-p_{1,1})_q!},\qquad  p_{2,2}\leq
p_{2,1}\,,
$$
$$
\Psi_{z_1,z_2}(p_{2,1},p_{2,2})=0,
\qquad p_{2,2}>p_{2,1}\,.
$$
Using the parametrization
$$
q=e^{-\epsilon},\qquad p_{21}=m(\epsilon)+x_{21}\epsilon^{-1}\qquad
p_{22}=-m(\epsilon)+x_{21}\epsilon^{-1}\qquad
z_i=e^{\imath\epsilon \lambda_i},\qquad   i=1,2,
$$
with $m(\epsilon)=-[\epsilon^{-1}\ln \,\epsilon]$ we obtain
$$
\Psi^{\mathfrak{gl}_2}_{z_1,z_2}(p_{21},p_{22})
=\sum_{x_{22}-\epsilon m(\epsilon)\leq x_{11}\leq
  x_{2,1}+\epsilon\,m(\epsilon)}\frac{e^{\imath\lambda_1 x_{11}+
\imath\lambda_2(x_{21}+x_{22}-x_{11})}}{((x_{11}-x_{22})/\epsilon+m(\epsilon))_q!\,
\,((x_{21}-x_{11})/\epsilon+m(\epsilon))_q!},
$$
where we use the notations $p_{11}=x_{11}/\epsilon$.
Taking into account
$$
\frac{1}{(y/\epsilon+m(\epsilon))_q!}
=e^{+\frac{\pi^2}{6}\frac{1}{\epsilon}+\frac{1}{2}\ln\frac{\e}{2\pi}
-e^{-y}+O(\epsilon)},
$$
we obtain
$$
\psi_{\lambda_1,\lambda_2}^{\mathfrak{gl}_2} (x_1,x_2)\,
  =\lim_{\epsilon\to +0}\,\,\,
  \epsilon\,e^{-\frac{\pi^2}{3}\frac{1}{\epsilon}-\ln\frac{\e}{2\pi}}\,\,
  \Psi^{\mathfrak{gl}_2}_{z_1,z_2}(p_{21},p_{22})\,
$$
$$
= \int_{\IR}\,dx_{11}e^{\imath\lambda_1 x_{11}}
e^{\imath\lambda_2(x_{21}+x_{22}-x_{11})}\,e^{-e^{x_{11}-x_{21}}-e^{x_{22}-x_{11}}}.
$$
\end{ex}

\vskip 1cm

\noindent {\small {\bf A.G.} {\sl Institute for Theoretical and
Experimental Physics, 117259, Moscow,  Russia; \hspace{8 cm}\,
\hphantom{xxx}  \hspace{2 mm} School of Mathematics, Trinity College
Dublin, Dublin 2, Ireland; \hspace{6 cm}\hspace{5 mm}\,
\hphantom{xxx}   \hspace{2 mm} Hamilton Mathematics Institute,
Trinity College Dublin, Dublin 2, Ireland;}}

\noindent{\small {\bf D.L.} {\sl
 Institute for Theoretical and Experimental Physics,
117259, Moscow, Russia};\\
\hphantom{xxxx} {\it E-mail address}: {\tt lebedev@itep.ru}}\\

\noindent{\small {\bf S.O.} {\sl
 Institute for Theoretical and Experimental Physics,
117259, Moscow, Russia};\\
\hphantom{xxxx} {\it E-mail address}: {\tt Sergey.Oblezin@itep.ru}}

\end{document}